\documentclass[12pt]{amsart}
\usepackage{amsmath,amscd,amssymb,amsfonts}
\def\ver{May 14, 2010, v.fin5}
\def\Mustata{Musta\c{t}\v{a}}
\def\1{\hskip1pt}
\def\bs{\bigskip}
\def\ms{\medskip}
\def\ssk{\smallskip}
\def\nin{\noindent}
\def\msum{\h{$\sum$}}
\def\mcap{\h{$\bigcap$}}
\def\mcup{\h{$\bigcup$}}

\def\mopl{\h{$\bigoplus$}}

\def\scirc{\,\raise.2ex\h{${\scriptstyle\circ}$}\,}
\def\ssb{\raise.2ex\h{${\scriptscriptstyle\bullet}$}}
\def\a{\alpha}
\def\b{\beta}
\def\A{{\mathcal A}}
\def\C{{\mathbf C}}

\def\dd{\langle d\rangle}
\def\e{\varepsilon}
\def\G{{\mathcal G}}
\def\h{\hbox}
\def\H{\widetilde{H}}
\def\I{{\mathcal I}}
\def\j{\tilde{j}}
\def\J{{\mathcal J}}
\def\LL{{\mathcal L}}
\def\la{\lambda}
\def\N{{\mathbf N}}
\def\OO{{\mathcal O}}
\def\P{{\mathbf P}}
\def\q{\quad}
\def\Q{{\mathbf Q}}
\def\tr{\tilde{\rho}}
\def\R{{\mathbf R}}
\def\S{{\mathcal S}}
\def\Si{\Sigma}
\def\Y{\widetilde{Y}}
\def\Z{\widetilde{Z}}
\def\ZZ{{\mathbf Z}}
\def\DR{{\rm DR}}
\def\IC{\h{\rm IC}}
\def\Gr{\h{\rm Gr}}
\def\Coker{\h{\rm Coker}}
\def\codim{{\rm codim}}
\def\mult{\h{\rm mult}}
\def\simto{\buildrel\sim\over\to}

\begin{document}
\title[First Milnor cohomology of hyperplane arrangements]
{First Milnor cohomology of hyperplane arrangements}
\author{Nero Budur}
\address{Department of Mathematics, The University of Notre Dame,
IN 46556, USA}
\email{nbudur@nd.edu}
\author{Alexandru Dimca}
\address{Laboratoire J.A.\ Dieudonn\'e, UMR du CNRS 6621,
Universit\'e de Nice-Sophia Antipolis, Parc Valrose,
06108 Nice Cedex 02, France}
\email{Alexandru.DIMCA@unice.fr}
\author{Morihiko Saito}
\address{RIMS Kyoto University, Kyoto 606-8502 Japan}
\email{msaito@kurims.kyoto-u.ac.jp}
\dedicatory{Dedicated to Professor Anatoly Libgober}
\date{\ver}
\begin{abstract}
We show a combinatorial formula for a lower bound of the dimension
of the non-unipotent monodromy part of the first Milnor cohomology
of a hyperplane arrangement satisfying some combinatorial conditions.
This gives exactly its dimension if a stronger combinatorial
condition is satisfied.
We also prove a non-combinatorial formula for the dimension of the
non-unipotent part of the first Milnor cohomology, which apparently
depends on the position of the singular points.
The latter generalizes a formula previously obtained by the second
named author.
\end{abstract}
\keywords{hyperplane arrangement, Milnor fiber, monodromy,
multiplier ideal}
\subjclass[2000]{32S22}
\thanks{The first author is supported by the NSF
grant DMS-0700360.}
\thanks{The second author is partially supported by
ANR-08-BLAN-0317-02 (SEDIGA)}
\thanks{The third author is partially supported by
Kakenhi 21540037.}
\maketitle

\centerline{\bf Introduction}

\bs
Let $D$ be a hyperplane arrangement in $X=\C^n$.
We assume $D$ is reduced and central, i.e.\ $D$ is defined by a
reduced homogeneous polynomial $f$.
We also assume that $D$ does not come from an arrangement in
$\C^{n-1}$, i.e.\ $f$ cannot be expressed as a polynomial of $n-1$
variables. 
Set $F_f=f^{-1}(1)$. This is the Milnor fiber of $f$.
Let $T$ be the Milnor monodromy on $H^j(F_f,\Q)$.
We have the monodromy decomposition
$$H^j(F_f,\C)=\mopl_{\la}\,H^j(F_f,\C)_{\la},$$
with $H^j(F_f,\C)_{\la}$ the $\la$-eigenspace.
Let $d=\deg D$. Set $Z:=\P(D)\subset Y:=\P^{n-1}$ and
$U:=\P^{n-1}\setminus Z$.
It is well-known (see e.g.\ [CS], [Di1]) that $H^j(F_f,\C)_{\la}=0$ if
$\la^d\ne 1$, and there are local systems $L^{(k)}$ of rank 1 on $U$
for $k=0,\dots,d-1$ such that
$$H^j(F_f,\C)_{\la}=H^j(U,L^{(k)})\q\h{for}\,\,\,
\la=\exp(2\pi ik/d),$$
and the monodromy around any irreducible component of $Z$ is the
multiplication by $\exp(-2\pi ik/d)$.
In particular, $L^{(0)}=\C_U$ so that $H^j(F_f,\C)_1=H^j(U,\C)$.
It has been conjectured that the $H^j(F_f,\C)_{\la}$ are combinatorial
invariants.
By [OS], the cohomology $H^{\ssb}(U,\C)=H^{\ssb}(F_f,\C)_1$ is
combinatorially described.

Let $Y'=\P^2$ be a sufficiently general linear subspace of dimension
2 in $Y=\P^{n-1}$.
By the generalized weak Lefschetz theorem for perverse sheaves
(i.e.\ Artin's theorem in [BBD]), we have an isomorphism
$$H^1(U,L^{(k)})=H^1(Y'\cap U,L^{(k)}|_{Y'\cap U}).\leqno(0.1)$$
Thus we may assume $n=3$ since we are interested in
$b_1(F_f)_{\la}:=\dim H^1(F_f)_{\la}$.

Set $\Si=\{y\in Z\mid\mult_yZ\ge 3\}$.
Let $Z_1,\dots,Z_d$ be the irreducible components of $Z$.
For $y\in\Si$ and $I\subset\{1,\dots,d\}$, set
$m_y:=\mult_yZ=\#\{i\mid y\in Z_i\}$, and

\ms\q\q\q\q\q
$m_{I,y}:=\#\{i\in I\mid y\in Z_i\},\q
\a_{I,y}:=\h{$\frac{|I|}{d}$}\1m_y-m_{I,y}$.

\ms
The following would be known to specialists.

\ssk\nin
{\bf Proposition~1.} {\it Let $\la=\exp(2\pi ik/d)$ with
$k\in\ZZ\cap[1,d/2]$.
Then $b_1(F_f,\C)_{\la}$ is described combinatorially by using
the Aomoto complex as in $(1.2)$ with $\omega$ given in $(1.2.3)$,
if there is
$I\subset\{1,\dots,d\}$ such that $|I|=k$ and

\ms\nin
$(0.2)$\q\q\q\q\q
either\q$\a_{I,y}\notin{\ZZ}_{>0}\,\,(\forall\,y\in\Si)$\q
or\q$\a_{I,y}\notin{\ZZ}_{<0}\,\,(\forall\,y\in\Si)$.}

\bs
In case the second condition of (0.2) is satisfied, we actually
consider $k':=d-k$ and the complement $I^c$ of $I$ instead of
$k$ and $I$ by using the complex conjugation on the Milnor
cohomology so that $\a_{I,y}$ is replaced by $-\a_{I,y}$.
In fact, this is the reason for which we assume $k\le d/2$;
the case $k>d/2$ is treated in the second case of (0.2).
The proof of Proposition~1 follows from a theorem of H.~Esnault,
V.~Schechtman and E.~Viehweg [ESV] about Aomoto's conjecture
(here [STV] is not needed since $n=3$) together with a combinatorial
description of the cohomology of the complement due to P.~Orlik and
L.~Solomon [OS], see (1.2) below.
Proposition~1 was essentially used in a calculation of examples
in [Sa2] although the combinatorial description of $H^{\ssb}(U)$ in
[OS] was not mentioned there.
It is possible to give a combinatorial description of $H^{\ssb}(U)$
in a more geometric way using the theory of perverse sheaves [BBD],
see also [BS].

If $k=d/3\in\N$ and $m_y=3\,(\forall\,y\in\Si)$, then
condition (0.2) in Proposition~1 becomes

\ms\nin
$(0.3)\,$\q\q\q\q\q
either\q$m_{I,y}\ge 1\,\,(\forall\,y\in\Si)$\q
or\q$m_{I,y}\le 1\,\,(\forall\,y\in\Si)$.

\ms\nin
Note that the first condition of (0.3) is equivalent to that $\Si$
is covered by $\bigcup_{i\in I}Z_i$.
There are examples where neither condition of (0.3) is satisfied,
see Example~(3.1)(iii) below.

\ssk
Returning to the general situation (with $n=3$), it is well-known
that $b_1(F_f)_{\la}=0$ unless there is $y\in\Si$ with $\la^{m_y}=1$,
see Remark~(3.4)(i) below.
(In this case we have actually a stronger assertion that
$b_1(F_f)_{\la}=0$ unless each $Z_i$ contains $y\in\Si$ with
$\la^{m_y}=1$, see Remark~(3.4)(ii) below.)
Using Proposition~1, we can show the following.

\ms\nin
{\bf Theorem~1.} {\it Assume $n=3$.
Let $m,r\in\ZZ\cap[3,d-1]$ with $d/m\in\ZZ$.
Assume there is a map $\phi:\{1,\dots,d\}\to\{1,\dots,r\}$.
Set $I_j=\phi^{-1}(j)$, $m_{j,y}=m_{I_j,y}$,
$\Si_{\dd}=\Si\setminus Z_d$,

\ms\q\q
$\Si^{\phi}=\{y\in\Si\mid \{y\}=Z_i\cap Z_{i'}\,\,\,\h{\it with}
\,\,\,\phi(i)\ne\phi(i')\},\q
\Si_{\dd}^{\phi}=\Si^{\phi}\cap\Si_{\dd}$.

\ms\nin
{\rm (i)} Assume $\la$ is an $m$-th root of unity and the following
conditions are satisfied.

\ssk\nin
$(0.4)\,\,\,$
$m_{j,y}>0\,\,(\forall\,y\in\Si_{\dd}^{\phi})$, and the
$m_{j,y}/m_{j',y}$ are independent of $y\in\Si_{\dd}^{\phi}$.

\ssk\nin
$(0.5)\,\,\,$
$\exists\,j_0\,$ s.t. $\,|I_{j_0}|=d/m$,
$\,\,m_{j_0,y}=\frac{1}{m}\1\msum_{j=1}^r\,m_{j,y}\,\,\,
(\forall\,y\in\Si_{\dd}^{\phi})$.

\ssk\nin
Then $b_1(F_f)_{\lambda}\ge r-2$.

\ms\nin
{\rm (ii)} Assume $\la$ is a primitive $m$-th root of unity and the
following condition is satisfied.

\ssk\nin
$(0.6)\,\,\,$
$m_{j,y}=1\,\,(\forall\,y\in\Si^{\phi})$, $m=r$, and $(0.2)$ is
satisfied for some $I_{j_0}$.

\ssk\nin
Then $b_1(F_f)_{\lambda}=m-2$.}

\ms
If the first condition of (0.6) is satisfied, i.e.\ if $m_{j,y}=1\,
(\forall\,y\in\Si^{\phi})$, then we have
$|\Si^{\phi}|=|I_j||I_{j'}|$ for any $j\ne j'$ so that $|I_j|=d/r$
for any $j$, and conditions (0.4-5) with $m=r$ are satisfied.
Hence $Z$ is an $(m,d/m)$-net in the sense of [FY].
By the Nullstellensatz there is a pencil such that
$\mcup_{\phi(i)=j}\,Z_i$ is a special fiber of the pencil for any
$j$, see [Yu2], Lemma~3.1.
This implies another proof of Theorem~1$\,$(i) in this case, see [DP],
Th.~3.1 (i) (or Remark~(3.3)(i) below).
Note that we have $m\le 4$ in this case by [St2], [Yu3].
For the moment any known examples are essentially of this type.
For $m=3$, there are lots of examples where (0.6) is satisfied,
see [Yu2] (and Remark~(3.2)(ii) below).
For $m=4$, however, only one example is known, see [FY] (and
Remark~(3.3)(iii) below).
Note that the last example implies a rather artificial example
where conditions (0.4--5) are satisfied for $m=4$ and $r=3$ by
considering $I_1$, $I_2$ and $I_3\cup I_4$ as a partition.
Note also that a hyperplane arrangement has a structure of a
multi-net in the sense of [FY] and $r\le 4$, if the hypotheses
(0.4--5) of Theorem~1(i) are satisfied, see Remark~(3.3)(ii) below.

\ms
In this paper we also prove a non-combinatorial formula for the
dimension of the non-unipotent monodromy part of the first Milnor
cohomology generalizing [Di1], Ch.6, Th.~4.15.
Let $\Si$ be as in Proposition~1, and set
$$\Si(k)=\{y\in\Si\mid m_yk/d\in\ZZ\}.$$
For $y\in\Si$, let $\I_{\{y\}}\subset\OO_Y$ be the reduced ideal of
$\{y\}\subset Y$, and define
$$\J^{(k)}:=\mcap_{y\in\Si}\,\I_{\{y\}}^{\lceil m_yk/d\rceil-2},\q
\J^{(>k)}:=\mcap_{y\in\Si}\,\I_{\{y\}}^{\lfloor m_yk/d\rfloor-1}.$$
Here $\lceil\a\rceil:=\min\{k\in\ZZ\mid k\ge\a\}$,
$\lfloor\a\rfloor:=\max\{k\in\ZZ\mid k\le\a\}$, and
$\I_{\{y\}}^j=\OO_Y$ for $j\le 0$.
Let $\C[X]_j$ denote the space of homogeneous polynomials
of degree $j$.
This is identified with $\Gamma(Y,\OO_Y(j))$. Define
$$J^{(k)}_j:=\Gamma(Y,\OO_Y(j)\otimes_{\OO_Y}\J^{(k)})\subset
\Gamma(Y,\OO_Y(j))=\C[X]_j.$$

\ms\nin
{\bf Theorem~2.} {\it Assume $n=3$.
For $k\in[1,d-1]$, let $k'=d-k$ and $\la=\exp(2\pi ik/d)$. Then
$$\aligned\dim\Gr_F^0H^1(F_f)_{\la}=\dim\Coker&\bigl(\rho^{(k)}:
J^{(k)}_{k-3}\to\mopl_{y\in\Si(k)}\,\J^{(k)}_y/\J^{(>k)}_y\bigr)\\
=\dim\Coker&\bigl(\tr^{(k)}:\C[X]_{k-3}\to\mopl_{y\in\Si}\,
\OO_{Y,y}/\J^{(>k)}_y\bigr),\\
\dim\Gr_F^1H^1(F_f)_{\la}=\dim\Coker&\bigl(\rho^{(k')}:J^{(k')}_{k'-3}
\to\mopl_{y\in\Si(k)}\,\J^{(k')}_y/\J^{(>k')}_y\bigr)\\
=\dim\Coker&\bigl(\tr^{(k')}:\C[X]_{k'-3}\to\mopl_{y\in\Si}\,
\OO_{Y,y}/\J^{(>k')}_y\bigr),\endaligned$$
and $b_1(F_f)_{\la}=\dim\Gr_F^1H^1(F_f)_{\la}+\dim\Gr_F^0H^1(F_f)_{\la}$.
Here we choose local trivializations of $\OO_Y(k-3)$, $\OO_Y(k'-3)$
to define the restriction morphisms $\tr^{(k)}$, $\tr^{(k')}$,
etc.\ at $y\in\Si$.}

\ms
Similar assertions in terms of superabundances (but without
reference to the mixed Hodge structure) were obtained by
A.~Libgober, see [Li1], [Li2].
Note that for $\a=\frac{k}{d}$ and $0<\e\ll 1/d$, we have
$$\J^{(k)}=\J(Y,(\a-\e)Z),\q\J^{(>k)}=\J(Y,\a Z),$$
where $\J(Y,\a Z)$ is the multiplier ideal sheaf [La],
see (2.3.4) below.
Moreover, the target of the restriction morphism $\rho^{(k)}$ is
identified with
$$\G(Y,\a Z):=\J(Y,(\a-\e)Z)/\J(Y,\a Z).$$

For simplicity assume $d/3\in\ZZ$ and $m_y=3\,(\forall\,y\in\Si)$.
Set $k=2d/3$, $k'=d/3$.
Then the target of $\rho^{(k')}$ vanishes, and $\rho^{(k)}$
coincides with $\tr^{(k)}$ which is identified with the evaluation
map (choosing points of $\C^3\setminus\{0\}$ representing the
points of $\Si$)
$$\mopl_{y\in\Si(k)}\,{\rm ev}^{k-3}_y:\C[X]_{k-3}\to
\mopl_{y\in\Si(k)}\,\C_y.$$
So we get a partial generalization of [Di1], Ch.6, Th.~4.15 where
$d=9$.
Note that ${\rm ev}^{k-3}_y$ can be defined and is surjective
if $2d/3<k<d$, although the surjectivity for $k=2d/3$ does not
hold in general.
This follows from formulas for the spectrum, see [BS], Th.~3 and 5
(these are closely related to {\Mustata}'s formula for the
multiplier ideals [Mu], see also [Te]).

\ms
We would like to thank A.~Libgober and S.~Yuzvinsky for useful
remarks, and the referee for valuable comments.

\ms
In Section 1 we recall theorems of Orlik, Solomon [OS] and
Esnault, Schechtman, Viehweg [ESV] to show Proposition~1, and
then prove Theorem~1.
In Section 2 we prove Theorem~2 using the theory of multiplier
ideals and Hodge theory.
In Section 3 we give some examples and remarks.

\bs\bs
\centerline{\bf 1. Combinatorial description}

\bs\nin
{\bf 1.1.~Orlik-Solomon algebra.}
Let $D$ be a central hyperplane arrangement in $X=\C^n$.
Set $Y=\P^{n-1}$, and $U=Y\setminus Z$, where
$Z=\P(D):=(D\setminus\{0\})/\C^*$.
Let $\S(Z)$ denote the intersection lattice consisting of any
intersections of the irreducible components $Z_i$ of $Z$.
For the definition of the Orlik-Solomon algebra, we allow the
ambient space $\P^{n-1}$ as a member, but the empty set
corresponding to $0\in\C^n$ is excluded.
(This is different from the definition in [BS], 1.1.)
In particular, $\S(Z)$ is not a lattice in the strict sense,
and is often called the ``poset of intersections" in the literature.
By [OS] there is an isomorphism of $\C$-algebras
$$A_{\S(Z)}^{\ssb}\simto H^{\ssb}(U,\C),$$
where $A_{\S(Z)}^{\ssb}$ is the quotient of the exterior algebra
$\bigwedge^{\ssb}\bigl(\bigoplus_i\C e_i\bigr)$ by the ideal
$\I_{\S(Z)}$.
Here the $e_i$ correspond to $Z_i$ for $i=1,\dots,d-1$ since we use
the induced affine arrangement on $\C^{n-1}=\P^{n-1}\setminus Z_d$,
see [Br].
Moreover, $\I_{\S(Z)}$ is determined by the combinatorial data,
see [OS].

\ms\nin
{\bf 1.2.~Solution of Aomoto's conjecture.}
Let $\a_i\in\C$ for $i\in\{1,\dots,d\}$, and assume
$\sum_{i=1}^d\a_i=0$, and
$$\msum_{Z_i\supset V}\,\a_i\notin{\bf N}\setminus\{0\}\q
\h{for any dense edge}\,\,\,V\in\S(Z).\leqno(1.2.1)$$
For the definition of ``dense", see [STV].
(In this paper the condition that $V$ is dense may be replaced with
$V\subset\Si$ since we assume essentially $n=3$.)
Note that (1.2.1) should be satisfied for $\S(Z)$, and not for
$\S(Z)_{\dd}$ in (1.3) below.
In the notation of (1.1), set
$$\omega=\msum_{i=1}^{d-1}\a_ie_i\in\A^1_{\S(Z)}.$$
This defines a complex $(\A^{\ssb}_{\S(Z)},\omega\wedge)$,
called the Aomoto complex associated to $\omega$.
Note that $e_i$ is identified with $dg_i/g_i$ where $g_i$ is a
linear function with a constant term defining
$Z_i\setminus Z_d\subset\C^{n-1}$.
We also have a regular singular connection $\nabla^{\omega}$
on $\OO_U$ such that
$$\nabla^{\omega}g=dg+g\omega\q\h{for}\,\,\,g\in\OO_U.$$
Note that (1.2.1) is the condition on the distribution of the
residues of the connection in [STV].
The main theorem of [ESV], [STV] assures that if (1.2.1) is satisfied,
then we have a quasi-isomorphism
$$H^{\ssb}_{\DR}(U,(\OO_U,\nabla^{\omega}))\simto
(\A^{\ssb}_{\S(Z)},\omega\wedge).\leqno(1.2.2)$$
This implies Proposition~1 in Introduction since the first condition
of (0.2) coincides with (1.2.1) by setting
$$\h{$\a_i=\frac{k}{d}-1\,$ if $\,i\in I$, and $\,\a_i=\frac{k}{d}\,$
otherwise.}\leqno(1.2.3)$$
In case the second condition of (0.2) is satisfied, we replace $k$
with $d-k$, and $I$ with its complement, using the complex
conjugation on the Milnor cohomology.

If condition (1.2.1) is not satisfied, then (1.2.2) may be false.
However, we have always the following inequality (see [LY], Prop.~4.2)
$$\dim H^j_{\DR}(U,(\OO_U,\nabla^{\omega}))\ge\dim
H^j(\A^{\ssb}_{\S(Z)},\omega\wedge).\leqno(1.2.4)$$
So we get for any $a\in\C^*$
$$\dim H^j_{\DR}(U,(\OO_U,\nabla^{a\omega}))\ge
\dim H^j(\A^{\ssb}_{\S(Z)},a\omega\wedge)=
\dim H^j(\A^{\ssb}_{\S(Z)},\omega\wedge).\leqno(1.2.5)$$

\ms\nin
{\bf 1.3.~Generalized residues.} Set
$$\S(Z)^{(k)}=\{V\in\S(Z)\mid\codim_YV=k\}.$$
Let $j_U:U\hookrightarrow Y$ denote the inclusion.
Since it is an affine morphism, $\R(j_U)_*\Q_U[n-1]$ is a perverse
sheaf [BBD] underlying naturally a mixed Hodge module [Sa1].
Let $W$ be the weight filtration.
By [BS], 1.7, there are constant variations of Hodge structures
$L_V$ of type $(0,0)$ on $V$ for $V\in\S(Z)$ such that we have
for $i=1,\dots,n-1$
$$\Gr^W_{n-1+i}(\R(j_U)_*\Q_U[n-1])=\mopl_{V\in\S(Z)^{(i)}}\,
L_V(-i)[n-1-i],\leqno(1.3.1)$$
where $\Gr^W_{n-1}(\R(j_U)_*\Q_U[n-1])=\Q_Y[n-1])$.
We may identify $L_V$ with a Hodge structure since it is constant.
Set
$$\S(Z)_{\dd}=\{V\in\S(Z)\mid V\not\subset Z_d\},\q
\S(Z)_{\dd}^{(k)}=\S(Z)_{\dd}\cap\S(Z)^{(k)}.$$
Then, restricting (1.3.1) to $Y\setminus Z_d=\C^{n-1}$, we get
by [BS], 1.9
$$H^i(U,\Q)=\mopl_{V\in\S(Z)_{\dd}^{(i)}}\,L_V(-i)\,\,\,
(i=1,\dots,n-1).\leqno(1.3.2)$$
This is compatible with the results of Brieskorn [Br] and
Orlik, Solomon [OS], since $\dim L_V$ is given by the M\"obius
function in loc.~cit.
For each $V\in \S(Z)_{\dd}^{(i)}$, we have the projection
$$\pi_V:H^i(U,\Q)\to L_V(-i).$$
This may be called the generalized residue along $V$.

In the notation of (1.2), set $e_i=dg_i/g_i$ for $1\le i<d$, and
$$e_{i_1,\dots,i_j}:=e_{i_1}\wedge\dots\wedge e_{i_j}\in
H^j(U,\Q).$$
This is compatible with the decomposition (1.3.2), i.e.
$$e_{i_1,\dots,i_j}\in L_V(-j)\subset H^j(U,\Q)\q\h{with}\,\,\,
V=Z_{i_1}\cap\cdots\cap Z_{i_j}.\leqno(1.3.3)$$
This means that its image by $\pi_{V'}$ for $V'\ne V$ vanishes.
The assertion is shown by taking a sub-arrangement $Z'\subset Z$
consisting of $Z_{i_1},\dots,Z_{i_j},Z_d$, and using the
compatibility of the exterior product with the restriction morphism
by the inclusion $U\hookrightarrow U':=\P^{n-1}\setminus Z'$,
since the construction of $L_V$ is functorial for $Z$.
Indeed, taking the graded pieces of the canonical morphism
$$\R(j_{U'})_*\Q_{U'}[n-1]\to\R(j_U)_*\Q_U[n-1],$$
the obtained morphism is compatible with the direct sum decomposition
in (1.3.1) (since the direct factors in (1.3.1) are simple perverse
sheaves with different supports).
Note that $H^j(U',\Q)=\Q$ where $j$ is as above, and $L_V=\Q$ if $Z$
is a divisor with normal crossings at the generic point of $V$.

\ms
The following lemma is well known to the specialists,
see [LY], Lemma~3.1 (and also [Fa], [Li2], [Yu1]) where the
situation is localized at $V$ using the fact that the relations
of the Orlik-Solomon algebra are of the form $\partial(e_J)$ for
certain $J$ and are compatible with the decomposition by $V$.
We note here a short proof for the convenience of the reader.

\ms
Set $I_V=\{i\mid Z_i\supset V\}$, and
$\a_V=\sum_{i\in I_V}\a_i$, where $\omega=\sum_{i=1}^{d-1}\a_ie_i$
as above.

\ms\nin
{\bf Lemma~1.4.} {\it
Let $\eta=\sum_{i=1}^{d-1}\b_ie_i$ with $\b_i\in\C$.
For $V\in\S(Z)_{\dd}^{(2)}$ we have}
$$\a_V\b_i=\b_V\a_i\,\,(\forall\,i\in I_V)\q\h{if}\,\,\,\,
\pi_V(\omega\wedge\eta)=0.\leqno(1.4.1)$$

\ms\nin
{\it Proof.}
For $i,j,k\in I_V$, it is easy to show the relation
(see e.g.\ [OS])
$$e_{i,j}+e_{j,k}=e_{i,k}.\leqno(1.4.2)$$
In the case $\codim_YV=2$, we have $\dim L_V=|I_V|-1$ by the
definition of $L_V$ in [BS], 1.7.
(This is related to the theory of M\"obius function [OS]).
Hence $L_V$ has a basis consisting of $e_{i,k}\,(i\ne k)$ for any
fixed $k\in I_V$.

By hypothesis we have for any $k\in I_V$ the vanishing of
$$\msum_{i,j\in I_V}\a_i\b_je_{i,j}
=\msum_{i,j\in I_V}\a_i\b_j(e_{i,k}-e_{j,k})
=\msum_{i\in I_V}(\a_i\b_V-\b_i\a_V)e_{i,k}.$$
This implies $\a_V\b_i=\b_V\a_i$ for any $i\in I_V\setminus\{k\}$.
However, $k$ can vary in $I_V$.
So the assertion follows.

\ms\nin
{\bf 1.5.~Kernel of the differential $\omega\wedge$.}
The subject of this subsection has been studied very well in [Fa],
[LY], [FY].
Assume $n=3$.
Set $J'=\{1,\dots,d-1\}$, and let
$$\omega=\msum_{i=1}^{d-1}\,\a_ie_i\q\h{with}\,\,\,\a_i\in\C^*.$$
With the notation of (1.4), set $\a_y=\a_V\,(:=\sum_{i\in I_V}\a_i)$
if $V=\{y\}$. Let $Z'=Z\setminus Z_d$, and
$$\Si_{\dd}=\{y\in Z'\mid\mult_yZ'\ge 3\},\q
\Si_{\dd}^{\a}=\{y\in\Si_{\dd}\mid\a_y=0\}.$$
Let $Z'{}^{\a}$ be the proper transform of $Z'$ by the blow-up
of $Y':=Y\setminus Z_d$ along $\Si_{\dd}^{\a}$.
We say that $I\subset J'$ is $\a$-connected, if $\mcup_{j\in I}\,
Z'_j$ is the image of a connected subvariety of $Z'{}^{\a}$,
where $Z'_1,\dots,Z'_{d-1}$ are the irreducible components of $Z'$.
Note that the set of $\a$-connected components is a
``neighborly partition" in [Fa], [LY].
Let $J'_1,\dots,J'_{r'}$ be the $\a$-connected components of $J'$.

Let $\eta=\sum_{i=1}^{d-1}\b_ie_i$ with $\b_i\in\C$, and assume
$$\omega\wedge\eta=0.$$
By Lemma~(1.4), there are $c_k\in\C\,\,(k=1,\dots,r')$ such that
$$\b_i=c_k\a_i\q\h{for}\,\,\,i\in J'_k.\leqno(1.5.1)$$
So we assume $r'\ge 2$ since $H^1(\A_{\S(Z)}^{\ssb},\omega\wedge)=0$
otherwise.
In this case there are intersections $Z'_i\cap Z'_{j}=\{y\}$
such that $i,j$ belong to different $\a$-connected components and
$y\in\Si_{\dd}^{\a}$.
So $\a_y=0$, and Lemma~(1.4) implies a further condition
$$\b_y:=\msum_{i\in I_y}\,\b_i=0\,\,(y\in\Si_{\dd}^{\a}).
\leqno(1.5.2)$$
However, the relation between these conditions for different
$y\in\Si_{\dd}^{\a}$ is quite nontrivial, see [Fa], [FY].
Since the image of $\omega\wedge$ in $\A_{\S(Z)}^1$ is
1-dimensional, we get at least
$$\dim H^1(\A_{\S(Z)}^{\ssb},\omega\wedge)\le r'-2.
\leqno(1.5.3)$$

\ms\nin
{\bf 1.6.~Proof of Theorem~1.}
Define the connection $\nabla^{\omega}$ as in (1.2) setting
$I=I_{j_0}$ for the $j_0$ in (0.5), where $k=d/m$ and the $\a_i$
are defined as in (1.2.3) (although (1.2.1) is not necessarily
satisfied).
Then $\a_y=0\,\,(y\in\Si_{\dd}^{\phi})$ by the last condition of (0.5),
and we get
$$\Si_{\dd}^{\a}=\Si_{\dd}^{\phi},$$
since the inclusion $\subset$ follows from the constantness of
$\a_i$ on each $I_j$.
Note that $I_j$ is not necessary $\a$-connected in the sense of
(1.5).
However, we may consider the case
$$\b_i=c_{\phi(i)}\,\,\,(i=1,\dots,d-1),$$
where $c_j\in\C\,\,(j=1,\dots,r)$.
Then the condition given by (1.5.2) is written as
$$\msum_{j=1}^r\,m_{j,y}c_j=0\,\,(y\in\Si_{\dd}^{\a}),$$
and it is independent of $y\in\Si_{\dd}^{\a}$ by the last condition
of (0.4). We get thus
$$\dim H^1(\A_{\S(Z)}^{\ssb},\omega\wedge)\ge r-2,$$
since $\dim\omega\wedge\A_{\S(Z)}^0=1$.
(Here $\Si_{\dd}^{\phi}$ may be empty.)
So the desired inequality for $\la=\exp(2\pi ia/m)\,\,(a=1,\dots,m)$
follows from (1.2.5).

If (0.6) is satisfied, then (0.4) and (0.5) are also satisfied.
In this case we have $\Si_{\dd}^{\phi}\ne\emptyset$ since
$|I_j|=d/r\ge 2$.
Since the Milnor monodromy is defined over $\Q$, it is enough to
consider the case $\la=\exp(2\pi i/m)$ using the Galois group of
$\overline{\Q}/\Q$ which acts transitively on the primitive $m$-th
roots of unity.
We define $\a_i$ and $\omega$ as above where (0.2) is satisfied for
$I=I_{j_0}$.
In this case each $I_j$ is $\a$-connected, since $Z'_i\cap Z'_{i'}$
for $i,i'\in I_j$ cannot belong to $\Si_{\dd}^{\a}$ by the condition
$m_{j,y}=1$ for $y\in\Si^{\phi}$.
So we get
$$\dim H^1(\A_{\S(Z)}^{\ssb},\omega\wedge)=m-2,$$
and the assertion follows from Proposition~1.
This finishes the proof of Theorem~1.

\bs\bs
\centerline{\bf 2. Non-combinatorial description}

\bs\nin
{\bf 2.1.~Mixed Hodge complexes.}
Set $\Si=\{y\in Z\mid\mult_yZ\ge 3\}$ as in Introduction.
Let $\pi:\Y\to Y$ be the blow-up of $Y=\P^2$ along $\Si$.
Set $\Z=\pi^{-1}(Z)$.
Let $L^{(k)}$ be the local system on $U$ calculating
$H^{\ssb}(F_f,\C)_{\la}$ for $\la=\exp(2\pi ik/d)$ as in Introduction.
Let $\LL^{(k)}$ be the Deligne extensions over $\Y$ such that the
eigenvalues of the residue of the connection are contained in
$[0,1)$, see [D1].

It is well-known (see e.g.\ [Es], [Ti]) that the Hodge
filtration $F$ on the Milnor cohomology $H^{\ssb}(F_f,\C)_{\la}=
H^{\ssb}(U,L^{(k)})$ is given by the logarithmic de Rham complex
whose $j$-th component is
$$\DR_{\log}^j\,\LL^{(k)}:=\Omega_{\Y}^j(\log\Z)
\otimes_{\OO}\LL^{(k)},$$
where the Hodge filtration $F^p$ is induced by the truncation
$\sigma_{\ge p}$ in [D2], 1.4.7.
By the strict compatibility of the Hodge filtration $F$, we get
$$\Gr_F^pH^{p+q}(F_f,\C)_{\la}=H^q(\Y,\Omega_{\Y}^p(\log\Z)
\otimes_{\OO}\LL^{(k)}).\leqno(2.1.1)$$
It is also well-known (see e.g.\ [BS], 1.4) that
$$\LL^{(k)}=\OO_{\Y}((-k)\H+\msum_{y\in\Si}\,\lfloor m_yk/d\rfloor
E_y),\leqno(2.1.2)$$
where $\lfloor\a\rfloor:=\max\{k\in\ZZ\mid k\le\a\}$, and $\H$
is the pull-back of a general hyperplane $H\subset Y$.

\ms\nin
{\bf 2.2.~Weight filtration.}
Let $\j:U\to\Y$ denote the inclusion.
The perverse sheaf $\R\j_*L^{(k)}[2]$ has the weight filtration $W$
such that
$$\Gr^W_2(\R\j_*L^{(k)}[2])=\IC_YL^{(k)}=(\j_*L^{(k)})[2],
\leqno(2.2.1)$$
where the middle term is the intersection complex with coefficients
in the local system $L^{(k)}$, see [BBD], [GM].
Here the local monodromy of $L^{(k)}$ is trivial only around the
exceptional divisors $E_y$ for $y\in\Si(k)$, and these are smooth
and disjoint.
This implies
$$\Gr^W_k(\R\j_*L^{(k)}[2])=0\q\h{for}\,\,\,k\ne 2,3.\leqno(2.2.2)$$
Let $\Z'$ be the proper transform of $Z$, and set
$E'_y=E_y\setminus\Z'$ with the inclusion $j_y:E'_y\to E_y$.
For $y\in\Si(k)$, $L^{(k)}$ can be extended over $U\cup E'_y$,
since the local monodromy around $E'_y$ is trivial.
Let $L^{(k)}_{E'_y}$ be its restriction over $E'_y$. Then
$$\Gr^W_3(\R\j_*L^{(k)}[2])=\mopl_{y\in\Si(k)}
((j_y)_*L^{(k)}_{E'_y})[1],\leqno(2.2.3)$$

\ms\nin
{\bf 2.3.~Multiplier ideals.}
Let $Y$ be a smooth complex algebraic variety, and $Z$ be
a divisor on it.
Let $\pi:\Y\to Y$ be an embedded resolution of $Z$.
Set $\Z:=\pi^*Z=\msum_i\,m_i\Z_i$.
The following is well-known (see e.g. [La]).

For $\a\in\Q_{>0}$, the multiplier ideal sheaves $\J(Y,\a Z)$
are defined by
$$\J(Y,\a Z)=\pi_*(\omega_{\Y/Y}\otimes_{\OO}\OO_{\Y}(-\msum_i\,
\lfloor\a m_i\rfloor \Z_i)),\leqno(2.3.1)$$
and we have the local vanishing theorem (see loc.~cit.\ 9.4.1)
$$R^i\pi_*(\omega_{\Y/Y}\otimes_{\OO}\OO_{\Y}(-\msum_i\,
\lfloor\a m_i\rfloor \Z_i))=0\,\,\,(i>0).\leqno(2.3.2)$$
Moreover, if $Y$ is proper and $Z'$ is another divisor on $Y$
such that $Z'-\a Z$ is nef and big, then we have the vanishing
theorem of Nadel (see loc.~cit.\ 9.4.8)
$$H^i(Y,\omega_Y\otimes_{\OO}\OO_{\Y}(Z')\otimes\J(Y,\a Z))=0\,\,\,
(i>0).\leqno(2.3.3)$$

If $Z$ is a projective hyperplane arrangement in $Y=\P^2$
with reduced structure, let $\I_{\{y\}}$ be the reduced ideal sheaf
of $\{y\}\subset Y$.
It is well-known that
$$\J(Y,\a Z)=\mcap_{y\in\Si}\,\I_{\{y\}}^{\lfloor\a m_y\rfloor-1}
\q\h{for}\,\,\,\a>0.\leqno(2.3.4)$$
This is a special case of {\Mustata}'s formula [Mu].

\ms\nin
{\bf Theorem~2.4.} {\it For $\la=\exp(2\pi ik/d)$, we have a
canonical isomorphism}
$$\Gr_F^0H^1(F_f,\C)_{\la}=H^1\bigl(Y,\OO_Y(k-3)\otimes_{\OO}\J
(\h{$\frac{k}{d}$}\1Z)\bigr)^{\vee}.\leqno(2.4.1)$$

\ms\nin
{\it Proof.}
By (2.1.1--2) together with Serre duality, we have
$$\aligned\Gr_F^0H^1(F_f,\C)_{\la}&=H^1\bigl(\Y,\OO_{\Y}((-k)\H+
\msum_{y\in\Si}\,\lfloor m_yk/d\rfloor E_y)\bigr)\\
&=H^1\bigl(\Y,\OO_{\Y}((k-3)\H+\msum_{y\in\Si}\,
(1-\lfloor m_yk/d\rfloor) E_y)\bigr)^{\vee}.\endaligned$$
Here $\omega_{\Y/Y}=\OO_{\Y}(\msum_{y\in\Si}E_y)$, and hence
$\Omega_{\Y}^2=\OO_{\Y}(-3\H+\msum_{y\in\Si}E_y)$.
So the assertion follows from (2.3.1--2).

\ms\nin
{\bf 2.5.~Proof of Theorem~2.}
Set first $\a=\frac{k}{d}$, and $\a'=\a-\e$ with $0<\e\ll 1/d$.
With the notation of Introduction, we have the short exact sequence
$$0\to\J(Y,\a Z)\to\J(Y,\a'Z)\to\G(Y,\a Z)\to 0.\leqno(2.3.1)$$
Since $\OO_Y(Z)=\OO_Y(d)$ and $\a'd<k$, we have by (2.3.3)
$$H^1(Y,\OO_Y(k-3)\otimes\J(Y,\a'Z))=0.\leqno(2.5.2)$$
Then, after taking the tensor product of (2.3.1) with $\OO_Y(k-3)$,
we get the associated long exact sequence
$$H^0(Y,\J(\a'Z)(k-3))\buildrel{\gamma}\over\to H^0(Y,\G(\a Z)
(k-3))\to H^1(Y,\J(\a Z)(k-3))\to 0,$$
where $\J(\a Z)(k-3):=\J(Y,\a Z)\otimes\OO_Y(k-3)$ and similarly for
$\G(\a Z)(k-3)$.
Moreover, $\gamma$ is identified with $\rho^{(k)}$ as is explained
after Theorem~2.
So the first equality in Theorem~2 follows from Theorem~(2.4).

For the proof of the second equality, we show, in the above
notation, that the first equality holds with the right-hand side
replaced by the dimension of
$$\Coker\bigl(H^0(Y,\J(\b Z)(k-3))\to H^0(Y,\J(Y,\b Z)/\J(Y,\a Z))
\bigr),$$
by decreasing induction on $0<\b<\a=\frac{k}{d}$.
If $\b=\a'$, this is the first equality.
Assume the assertion holds for $\b$, and set $\b'=\b-\e$.
We have a long exact sequence
$$0\to H^0(Y,\J(\b Z)(k-3))\to H^0(Y,\J(\b'Z)(k-3))
\buildrel{\gamma}\over\to H^0(Y,\G(\b Z)(k-3))\to,$$
and the morphism $\gamma$ is surjective by (2.5.2) with $\a'$
replaced by $\b$.
So the assertion holds also for $\b'$ using the snake lemma
since the support of $\J(Y,\b'Z)/\J(Y,\a Z)$ is contained in a
0-dimensional subset $\Si$.
Thus the assertion for $\Gr_F^0H^1(F_f)_{\la}$ is proved.

For $\Gr_F^1H^1(F_f)_{\la}$, we have by (2.2.2--3)
$$\h{$H^1(F_f)_{\la}\oplus H^1(F_f)_{\overline{\la}}\,$ is pure of
weight $1$.}\leqno(2.5.3)$$
Indeed, using the weight spectral sequence, it is enough to show
that the local system $L^{(k)}_{E'_y}$ on $E'_y$ has no nontrivial
global sections, i.e.\ the local monodromies of $L^{(k)}_{E'_y}$
around $E_y\cap\Z'$ are nontrivial.
But this is clear by the construction of $L^{(k)}_{E'_y}$ before
(2.2.3).
So (2.5.3) follows.
(See [DP] for another proof of (2.5.3).)

Using the complex conjugation, we then get
$$\dim\Gr_F^1H^1(F_f)_{\la}=\dim\Gr_F^0H^1(F_f)_{\overline{\la}},$$
and the assertion for $\Gr_F^1H^1(F_f)_{\la}$ follows.
This finishes the proof of Theorem~2.

\bs\bs
\centerline{\bf 3. Examples and Remarks}

\bs\nin
{\bf Examples~3.1.}
The following examples are studied in [CS], [Di1], etc.

\ms\nin
(i) Let $n=3$, $d=6$, and
$$f=xyz(x-y)(x-z)(y-z).$$
This is the simplest example with $b_1(F_f)_{\la}\ne 0\,\,
(\la\ne 1)$.
Here $\la=\exp(2\pi i/3)$.
In this case, (0.2) for $k=2$ and (0.6) for $r=m=3$ are both
satisfied.

\ms\nin
(ii) Let $n=3$, $d=9$, and
$$f=xyz(x-y)(y-z)(x-y-z)(2x+y+z)(2x+y-z)(2x-5y+z).$$
This is the second simplest example with $b_1(F_f)_{\la}\ne 0\,
(\la\ne 1)$.
Here $\la=\exp(2\pi i/3)$.
In this case, (0.2) for $k=3$ and (0.6) for $r=m=3$ are both
satisfied.
This arrangement is dual to the Pappus configuration.

\ms\nin
(iii) Let $n=3$, $d=9$, and
$$f=xyz(x+y)(y+z)(x+3z)(x+2y+z)(x+2y+3z)(2x+3y+3z).$$
Then $b_1(F_f)_{\la}=0$ for $\la\ne 1$.
In this case neither (0.2) for $k=3$ nor (0.6) for $r=m=3$ is
satisfied.
So we need Theorem~2 or [Di1], Ch.6, Th.~4.15 to calculate
$b_1(F_f)_{\la}$ for $\la=\exp(\pm 2\pi i/3)$, see also [CS].
The referee has pointed out that (1.2.2) holds for this example
even though (0.2) is not satisfied.

\ms\bs\nin
{\bf Remarks~3.2.}
(i) According to A.~Libgober, there is an example with $n=3$,
$d=9$, $|\Si|=12$, and $Z:=\P(D)$ has no ordinary double point.
It is the dual of the nine inflection points of a smooth cubic
curve $E$ in the dual projective space $\P^2$, see [Li2], p.~243,
Ex.~2.
This is shown by choosing an inflection point as the origin $O$
of the group law so that there is a line passing through three
points $P_1,P_2,P_3$ on $E$ if and only if
$$P_1+P_2+P_3=O.$$
In this case conditions (0.4--5) are satisfied so that
$b_1(F_f)_{\la}\ne 0$ for $\la=\exp(\pm 2\pi i/3)$ although
(0.2) cannot to be satisfied.
The referee has pointed out that the inequality (1.2.4) is strict in
this example.

If the cubic curve is given by $u^3+v^3+w^3=0$, then
$$f=(x^3-y^3)(x^3-z^3)(y^3-z^3),$$
and setting $\theta=\exp(\pm 2\pi i/3)$, we have
$$\Si=\{(\theta^i:\theta^j:1)\mid i,j=0,1,2\}\cup
\{(1:0:0),\,(0:1:0),\,(0:0:1)\}\subset\P^2.$$
In this case the evaluation map $\mopl_{y\in\Si}\,{\rm ev}_y^3$
explained after Theorem~2 is injective so that $b_1(F_f)_{\la}=2$
for $\la=\exp(\pm 2\pi i/3)$ using [Di1], Ch.6, Th.~4.15 or Theorem~2
in Introduction.
(This follows from the assertion that $a=b=c=0$ if
$a\theta^{2i}+b\theta^i+c=0$ for $i=0,1,2$ where $a,b,c\in\C$.
We have a more geometric proof using the fact that any $Z_i$
contains 4 points of $\Si$ so that any cubic polynomial vanishing
on $\Si$ has to vanish on any $Z_i$.)

\ms
(ii) The argument in Remark~(i) can be used to construct further
examples with $b_1(F_f)_{\la}=1$ for $\la=\exp(\pm 2\pi i/3)$.
Indeed, take a finite subgroup $G$ of $E$ together with three points
$P_1,P_2,P_3$ such that $P_1+P_2+P_3=O$ and $P_i-P_j\notin G$ for
any $i\ne j$, see [Yu2].
Consider the dual line arrangement associated to $P_i+Q$ for
$i=1,2,3$ and $Q\in G$ where $d=3|G|$ and $\mult_yZ\le 3$ for any
$y\in\Si$ (since the degree of the curve is 3).
Then (0.6) for $r=m=3$ is satisfied if furthermore $3P_i\notin G$
for some $i$ (since the last condition implies that the second
condition of (0.3) holds).
This gives examples with the same combinatorial data as
Examples~(3.1)(i) and (ii).
Here ``same combinatorial data" means that there is an isomorphism
between the intersection lattices defined by the intersections of
their irreducible components.
If we set $P_1=0$ with $P_2=-P_3$ generic and if $|G|=3$, then we
get an example where $d=9$, $|\Si|=10$ and $b_1(F_f)_{\la}=1$ for
$\la=\exp(\pm 2\pi i/3)$ by Theorem~1$\,$(ii).
This example is the specialization of the Pappus arrangement in
Example~(3.1)(ii), see [Fa].

\ms
(iii) As for Example~(3.1)(iii), a line arrangement with the same
combinatorial data can be constructed as follows.
Choose a finite subgroup $G=\ZZ/27\ZZ$ in an elliptic curve, and
then take the subset $J$ of $G$ consisting of $a\in G$ with
$a-1\in 3G$.
The irreducible components of $D$ defined by the linear factors
of $f$ in Example~(3.1)(iii) correspond, for example, to
$$7,\,1,\,4,\,19,\,22,\,16,\,13,\,10,\,25\in\ZZ/27\ZZ,$$
respecting the order of the factors of $f$
(and this is confirmed by the referee using Mathematica).
It means that we have $a+b+c=0$ in $G$ if and only if the
corresponding three lines meet at a point.
(We can associate $a':=(a-1)/3\in\ZZ/9\ZZ$ to $a\in J$.
Then the above condition is equivalent to that $a'+b'+c'=-1$ in
$\ZZ/9\ZZ$, and the latter may be easier to handle.)

\ms
(iv)
In the case $\mult_yZ\le 3\,\,(\forall\,y\in\Si)$, it is sometimes
possible to construct a finite subset of an elliptic curve defining
a line arrangement which is combinatorially equivalent to a given
line arrangement $Z$ as follows.
Let $M$ be a matrix of size $(q,d)$ (where $q=|\Si|$) such that
if $y\in\Si$ is the intersection of three lines $Z_{i_1}$,
$Z_{i_2}$, $Z_{i_3}$, then it corresponds to a row vector $v_y$ of
$M$ such that
$$v_{y,i}=\begin{cases}1&\h{if}\,\,\,i\in\{i_1,i_2,i_3\},\\
0&\h{if}\,\,\,i\notin\{i_1,i_2,i_3\}.\end{cases}$$
Consider the equation
$$Mx=0.$$
If it has a nontrivial solution $x$ in a finite
abelian group $G$ generated by two elements,
and if the components $x_i$ of $x$ are all different from each other,
then the $x_i$ may define the desired line arrangement after
embedding $G$ into an elliptic curve.
Here one problem is that new triple points may appear,
i.e.\ there may exist distinct elements $x_{i_1},x_{i_2},x_{i_3}$
whose sum is $O$ and which do not correspond to any row vector of
$M$.
(In the case of Example~(iii) the size of $M$ is $9\times 9$, and
its determinant is 27.
So we can solve the equation mod 27.
We can also count the number of the triple points, which is 9.)

\ms\nin
{\bf Remarks~3.3.}
(i) If $m_{j,y}=1\,(\forall\,y\in\Si^{\phi})$ in the notation of
Theorem~1, then $|I_j|$ is independent of $j$ as remarked after
Theorem~1, and hence $Z$ is an $(m,d/m)$-net.
So it is associated to a pencil on $\P^2$ by the Nullstellensatz,
see [Yu2], Lemma~3.1.
Then Theorem~1$\,$(i) follows as is shown in [DP], Th.~3.1 (i) in a
more general case.
Indeed, the existence of a pencil implies that there is a
surjective morphism of algebraic varieties
$$h:U\to S:=\P^1\setminus\{\h{$m$ points}\},$$
and moreover the local system $L^{(k)}$ on $U$ is the pull-back
of a local system $L_S^{(k)}$ on $S$.
Since $\dim H^1(S,L_S^{(k)})=m-2$, Theorem~1$\,$(i) then follows
from the injectivity of the pull-back
$$H^1(S,L_S^{(k)})\to H^1(U,L^{(k)}).$$
Note that the condition $m_{j,y}=1\,(\forall\,y\in\Si^{\phi})$
implies that there is a natural compactification
$\bar{h}:\bar{U}\to S$ of $h:U\to S$ such that
$\bar{U}\setminus U=\Si^{\phi}\times S$.
Hence the base change holds for $h:U\to S$ and
$\{s\}\hookrightarrow S$, and the Leray spectral sequence for $h$
degenerates at $E_2$ since $E_2^{p,q}=0$ unless $q=0,1$.

\ms
(ii) By [DPS], [FY], [Yu3], we have $r\le 4$ if conditions (0.4--5)
are satisfied in Theorem~1$\,$(i).
Indeed, if $E$ denotes the irreducible component of the resonance
variety $R_1(Z)$ (see [Fa]) containing $\omega$ in the proof of
Theorem~1, then we have by [DPS], Prop.~7.8(i)
$$\dim E=H^1(\A_{\S(D)}^{\ssb},\omega\wedge)+1\ge r-1,$$
where the last inequality follows from the proof of Theorem~1(i).
Moreover, there is a correspondence between the components $E$ and
the pencils, see [FY] (and also [Di3]).
Then we get the inequality $\dim E\le 3$ by using [Yu3].
So the desired inequality follows.
This argument also shows that a hyperplane arrangement has a
structure of a multi-net in the sense of [FY] if the hypotheses
(0.4) and (0.5) of Theorem~1(i) are satisfied.

\ms
(iii) For the moment only one example is known where the hypotheses
of Theorem~1 are satisfied with $m=4$.
This is induced by the Hesse pencil (see [FY], Ex.~3.5 and [St1])
$$s_0(x^3+y^3+z^3)-s_1xyz\q(s_0,s_1\in\C).$$
Calculating the logarithmic differential of $(x^3+y^3+1)/xy$,
the singular members of the pencil are given by
$$s_1/s_0=3,\,\,3\theta,\,\,3\bar{\theta},\,\,\infty,$$
where $\theta=\exp(\pm 2\pi i/3)$.
So the line arrangement is defined by
$$\h{$f=xyz\prod_{j=0}^2(x^3+y^3+z^3-3\theta^jxyz)=
xyz\prod_{i,j=0}^2(\theta^ix+\theta^jy+z)$}.$$
In this case $d=12$, and $|\Si|=9$. More precisely
$$\aligned\Si&=\{x^3+y^3+z^3=0\}\cap\{xyz=0\}\\
&=\{(-\theta^i:1:0),\,(-\theta^i:0:1),\,(0:-\theta^i:1)\mid
i=0,1,2\}\subset\P^2.\endaligned$$
Then condition (0.6) is satisfied for $m=r=4$ where
$\phi$ corresponds to the first factorization of $f$.
So we get $b_1(F_f)_{-1}\ge 2$ by Theorem~1$\,$(i),
and $b_1(F_f)_{\pm i}=2$ by Theorem~1$\,$(ii) (see also [DP]).
For $\la=-1$, we can calculate $b_1(F_f)_{-1}$ by using Theorem~2
where $d=12$, $k=k'=6$, and $m_y=4$ for $y\in\Si$.
In this case, $\dim\C[X]_3=10$, and $\rho^{(6)}=\tr^{(6)}$ is
identified with the evaluation map at $\Si$.
Its kernel is generated by $x^3+y^3+z^3$ and $xyz$, and its
cokernel is 1-dimensional.
So we get $b_1(F_f)_{-1}=2$.
It is also possible to calculate $b_1(F_f)_{\pm i}$ by using
Theorem~2.

\ms\nin
{\bf Remarks~3.4.} (i)
Let $f$ be a holomorphic function on a complex manifold $X$.
Let $\psi_{f,\la}\C_X$ denote the $\la$-eigenspace of the nearby
cycle functor $\psi_f\C_X$.
It is well-known that $\psi_{f,\la}\C_X$ is a perverse sheaf up to
a shift of complex, and its stalk at $y\in f^{-1}(0)$ gives
$H^{\ssb}(F_{f,y},\C)_{\la}$ where $F_{f,y}$ denotes the Milnor
fiber of $f$ around $y$.
These imply that $H^j(F_{f,y},\C)_{\la}=0$ for $j\ne\dim X-1$ if
$H^j(F_{f,y'},\C)_{\la}=0$ for any $j\in\ZZ$ and any $y'\ne y$
sufficiently near $y$.

\ms
(ii) In the case of a hyperplane arrangement, an assertion stronger
than (i) is known as follows.
We have $b_1(F_f)_{\la}=0$ unless each $Z_i$ contains $y\in\Sigma$
with $\la^{m_y}=1$ (see [Li3], Th.~3.1 or [Di2], Th.~6.4.13),
where we may assume $n=3$ by (0.1).
Indeed, if there is $Z_i$ containing no point $y$ of $\Si$ with
$\la^{m_y}=1$, then letting $j:U\hookrightarrow Y$,
$j_i:U_i:=Y\setminus Z_i\hookrightarrow Y$ denote the inclusions,
we have
$$\R j_*L^{(k)}=(j_i)_!j_i^*(\R j_*L^{(k)}),$$
using the blow-up along $\Si$, where $\la=\exp(2\pi ik/d)$.
Then the assertion follows from Artin's theorem on the
generalization of the weak Lefschetz theorem for perverse sheaves
[BBD].

\end{document}